\documentclass[12pt,a4paper,twoside]{article}
\usepackage[latin1]{inputenc}
\usepackage{amsmath,amsfonts,amssymb,a4,enumerate,epsfig,array,theorem}

\def\qed{\unskip\nobreak\hfil\penalty50\hskip1.75em\null\nobreak\hfil
$\blacksquare$ {\parfillskip=0pt \finalhyphendemerits=0 \par}\goodbreak}

\newfont{\eightrm}{cmr8}
\newfont{\fiverm}{cmr5}

\newcommand\proof{{\noindent\bf Proof: }}
\newcommand{\card}{\operatorname{card}}
\newcommand{\trace}{\operatorname{trace}}
\newcommand{\Ab}{{\operatorname{Ab}}}

\newcommand{\TT}{{\mathbb{T}}}
\newcommand{\ZZ}{{\mathbb{Z}}}
\newcommand{\QQ}{{\mathbb{Q}}}
\newcommand{\RR}{{\mathbb{R}}}
\newcommand{\CC}{{\mathbb{C}}}
\newcommand{\Ss}{{\mathbb{S}}}
\newcommand{\Ee}{{\cal E}}
\newcommand{\Ff}{{\cal F}}
\newcommand{\Gg}{{\cal G}}
\newcommand{\Jj}{{\cal J}}
\newcommand{\Ll}{{\cal L}}
\newcommand{\zz}{{\bf z}}
\newcommand{\ww}{{\bf w}}
\newcommand{\bT}{{\bf T}}
\newcommand{\bF}{{\bf F}}

\newtheorem{theo}{Theorem}
\newtheorem{prop}{Proposition}[section]
\newtheorem{coro}[prop]{Corollary}

\newtheorem{lemma}[prop]{Lemma}

\begin{document}
\title{Stability of compact actions \\ of the Heisenberg group}
\author{ Tania M. Begazo and Nicolau C. Saldanha }
\maketitle

\begin{abstract}
Let $G$ be the Heisenberg group of real lower triangular $3 \times 3$
matrices with unit diagonal.
A locally free smooth action of $G$ on a manifold $M^4$ is
given by linearly independent vector fields $X_1, X_2, X_3$ such that
$X_3 = [X_1,X_2]$ and $[X_1,X_3] = [X_2,X_3] = 0$.
The $C^1$ topology for vector fields induces
a topology in the space of actions of $G$ on $M^4$.
An action is {\sl compact} if all orbits are compact.
Given a compact action $\theta$,
we investigate under which conditions its $C^1$ perturbations $\tilde\theta$
are guaranteed to be compact.
There is more than one interesting definition of stability,
and we show that in the case of the Heisenberg group,
unlike for actions of $\RR^n$,
the definitions do not turn out to be equivalent.
\end{abstract}

{\noindent\bf Keywords:} Heisenberg group, nilpotent, compact action, stability.

{\noindent\bf MSC-class:} 37C85; 57R30; 57M60

\section{Introduction}
Let $G$ be the Heisenberg group of real lower triangular $3 \times 3$
matrices with unit diagonal:
$G$ is a nonabelian nilpotent Lie group of dimension $3$.
Let $\Gg$ be its Lie algebra of strictly lower triangular matrices.
Let $M_0$ be a smooth manifold of dimension $4$:
a locally free smooth action $\theta: G \times M_0 \to M_0$ is
given by linearly independent vector fields $X_1, X_2, X_3$ such that
$X_3 = [X_1,X_2]$ and $[X_1,X_3] = [X_2,X_3] = 0$.
A Riemann metric on $M_0$ induces a $C^1$ metric for vector fields on $M_0$
and therefore a $C^1$ metric in the space of such actions.

An action $\theta$ is {\sl compact} if all its orbits are compact.
Is it true that all actions $\tilde\theta$ which are
near $\theta$ (in the $C^1$ metric) are also compact?
This is true for some, but not all, compact actions $\theta$ 
and the best way to attack this problem is to consider 
first a neighborhood of an orbit and actions of a very special form. 

A given homomorphism $\phi: G \to G$ may be described as follows.
Take the 1-parameter subgroups $\alpha_1, \alpha_2, \alpha_3: \RR \to G$,
\[ \alpha_1(t) = \phi\begin{pmatrix}1 & 0 & 0 \\
t & 1 & 0 \\ 0 & 0 & 1\end{pmatrix},\quad
\alpha_2(t) = \phi\begin{pmatrix}1 & 0 & 0 \\
0 & 1 & 0 \\ 0 & t & 1\end{pmatrix},\quad
\alpha_3(t) = \phi\begin{pmatrix}1 & 0 & 0 \\
0 & 1 & 0 \\ t & 0 & 1\end{pmatrix}. \]
Let
\[ \alpha'_j(0) = \begin{pmatrix} 0 & 0 & 0 \\
a_{1j} & 0 & 0 \\ a_{3j} & a_{2j} & 0 \end{pmatrix}; \]
the numbers $a_{ij}$ are entries of a $3 \times 3$ matrix $A$
with $a_{13} = a_{23} = 0$, $a_{33} = a_{11} a_{22} - a_{12} a_{21}$.
Conversely, a matrix with these properties defines a homomorphism.
The homomorphism $\phi$ is an automorphism exactly when $A$ is invertible.

Let $H \subset G$ be the group of matrices in $G$ with integer coordinates:
$H$ is a discrete cocompact subgroup of $G$.
Let $M = (G/H) \times (-\epsilon,\epsilon)$
where $G/H = \{gH, g \in G\}$ is a compact manifold of dimension $3$.
An action is {\em homogeneous horizontal} if it is of the form
\[\begin{matrix}
\theta: & G \times M & \to & M \\
& (g, (g_1H, z)) & \mapsto & (\phi_z(g) g_1H, z)
\end{matrix}\]
where $\phi_z$ is a smooth family of automorphisms of $G$
such that $\phi_0$ is the identity.
The orbits of $\theta$ are of the form $(G/H) \times \{z\}$
and therefore compact.
With the notation of the previous paragraph, $\theta$ is described
by a family of matrices $A(z)$, $z \in (-\epsilon, \epsilon)$
satisfying $A(0) = I$,
$a_{13}(z) = a_{23}(z) = 0$,
$a_{33}(z) = a_{11}(z) a_{22}(z) - a_{12}(z) a_{21}(z)$.
Let
\[ A^\sharp(z) = \begin{pmatrix} a'_{11}(z) & a'_{12}(z) \\
a'_{21}(z) & a'_{22}(z) \end{pmatrix}. \]

The orbit $(G/H) \times \{0\}$ of a homogeneous horizontal action
$\theta$ is ($C^1$) {\em totally stable} (or {\em T-stable})
if for any neighborhood $U$ of $(G/H) \times \{0\}$ there is $\delta > 0$
such that if $\tilde X_1, \tilde X_2, \tilde X_3$ are smooth vector fields
on $U$ satisfying 
\[ \tilde X_3 = [\tilde X_1, \tilde X_2], \quad
[\tilde X_1, \tilde X_3] = [\tilde X_2, \tilde X_3] = 0, \quad
d_{C^1}(X_j, \tilde X_j) < \delta, (j = 1,2,3) \]
then the orbit by $\tilde X_1, \tilde X_2, \tilde X_3$
of any point of the form $(\ast,0)$ is compact.
Otherwise, the orbit is {\em T-unstable}.

The orbit $(G/H) \times \{0\}$ of a homogeneous horizontal action
$\theta$ is ($C^1$) {\em locally unstable} (or {\em L-unstable})
if for any $\delta > 0$ there exist smooth vector fields 
$\tilde X_1, \tilde X_2, \tilde X_3$ in $(G/H) \times (-\epsilon,\epsilon)$
satisfying 
\[ \tilde X_3 = [\tilde X_1, \tilde X_2], \quad
[\tilde X_1, \tilde X_3] = [\tilde X_2, \tilde X_3] = 0, \quad
d_{C^1}(X_j, \tilde X_j) < \delta, (j = 1,2,3) \]
such that $\tilde X_j$ coincides with $X_j$ outside
$(G/H) \times (-\delta,\delta)$
and such that the orbits by $\tilde X_1, \tilde X_2, \tilde X_3$
of any point of the form $(\ast,0)$ is noncompact.
Otherwise, the orbit is {\em L-stable}.

Our main theorems concern the stability of
the orbit $(G/H) \times \{0\}$ of a homogeneous horizontal action $\theta$.

\begin{theo}
\label{theo:L-unstable}
If $A^\sharp(0)$ is nilpotent then the orbit is L-unstable.
\end{theo}

\begin{theo}
\label{theo:mixed}
If \[ A'(0) = \begin{pmatrix} -\lambda & 0 & 0 \\ 0 & 2\lambda & 0 \\
0 & 0 & \lambda \end{pmatrix}, \lambda \ne 0 \]
then the orbit is L-stable and T-unstable.
\end{theo}

\begin{theo}
\label{theo:L-stable}
If $\det A^\sharp(0) \ne 0$ then the orbit is L-stable.
\end{theo}

\begin{theo}
\label{theo:T-stable}
If $A^\sharp(0)$ has no real eigenvalues then the orbit is T-stable.
\end{theo}

In section 2 we justify our interest in homogeneous horizontal actions
by proving that for any compact action, 
the neighborhood of an orbit admits a finite covering
which is homogeneous horizontal.
In section 3 we prove the instability claims, i.e., 
we present explicit perturbations of certain actions.
The main result of section 4 is that a certain vector field
is uniquely ergodic.
In section 5 we study the underlying foliations.
Finally, in section 6, we prove our main results.

This paper follows closely the first author's thesis (\cite{B}),
except for a chapter which was developed separately (\cite{BS}).
The simpler case of actions of $\RR^n$ was studied earlier by the
second author (\cite{S}).

The authors acknowledge support from CNPq, CAPES and Faperj (Brazil).
This work was done while the first author was in a leave of absence
from Universidade Federal do Par\'a and had a visiting position in UFF.

\section{Homogeneous horizontal actions}

Let $G$ be the Heisenberg group,
$H \subset G$ be the discrete subgroup of integer matrices,
$M_0$ be an arbitrary $4$-manifold and $M = G/H \times (-\epsilon, \epsilon)$.
\begin{prop}
\label{prop:homohoriz}
Let $\theta_0: G \times M_0 \to M_0$ be a smooth locally free compact action.
For any $p_0 \in M_0$ there exists an open neighborhood $V$ of $p_0$
with the following properties:
\begin{enumerate}
\item{$V$ is invariant under $\theta_0$;}
\item{there exists a covering map of finite order $\psi: M \to V$
and a homogeneous horizontal action $\theta$ on $M$ such that
$\theta_0(g,\psi(p)) = \psi(\theta(g,p))$ for all $g \in G$ and $p \in M$.}
\end{enumerate}
\end{prop}

\proof
Since orbits are compact and have codimension $1$,
the holonomy of each orbit of $\theta_0$ contains
at most two elements and admits a saturated open neighborhood
(see \cite{CL} or \cite{Camacho}).
Take $V$ to be such a neighborhood.
We assume the holonomy to be trivial;
otherwise just take a double cover.

Let $\gamma: (-\epsilon, \epsilon) \to M_0$ be a smooth curve
transversal to $\theta$ with $\gamma(0) = p_0$.
For each $z$, the subgroup $H_z \subset G$ of isotropy of $\gamma(z)$
is cocompact and discrete:
let $A_z = H_z/[G,G] \subset G/[G,G] \equiv \ZZ^2$ and
$a(z), b(z) \in H_z$ be continuously chosen elements projecting
via the quotient by $[G,G]$ to generators of $A_z$.
Let $\phi_z: G \to G$ be the automorphism with 
\[ \phi_z\left( \begin{pmatrix} 1 & 0 & 0 \\ 1 & 1 & 0 \\ 0 & 0 & 1
\end{pmatrix}\right) = a(z), \qquad
\phi_z\left( \begin{pmatrix} 1 & 0 & 0 \\ 0 & 1 & 0 \\ 0 & 1 & 1
\end{pmatrix}\right) = b(z) \]
and define $\psi$ and $\theta$ by
$\psi((H,z)) = \gamma(z)$ and 
$\theta(g_1, (g_2H, z)) = (\phi_z^{-1}(g_1) g_2H, z)$
so that $\psi((gH,z)) = \theta_0(\phi_z(g), \gamma(z))$. 
We are left with proving that $\psi$ is indeed a covering map of finite order
and that $\theta$ is well defined: these are rather straightforward checks
and are left to the reader.
Notice that the order of $\psi$ as a covering is the index in $H_z$
of the subgroup generated by $a(z), b(z)$.
\qed

\section{Instability}

Points of $G/H$ have coordinates $(y_1, y_2, y_3)$, meaning
\[ \begin{pmatrix} 1 & 0 & 0 \\ y_1 & 1 & 0 \\ y_3 & y_2 & 1 \end{pmatrix} H \]
so that $(y_1, y_2, y_3) = (y'_1, y'_2, y'_3)$ if and only if
$y'_1 - y_1$, $y'_2 - y_2$ and $(y'_3 - y_3) - y_2(y'_1 - y_1)$ are integers.
Alternatively, $G/H$ can be constructed from the cube $[0,1]^3$
(with points $(y_1, y_2, y_3)$) by glueing faces as follows:
$(y_1,y_2,0)$ with $(y_1,y_2,1)$; $(y_1,0,y_3)$ with $(y_1,1,y_3)$;
$(0,y_2,y_3)$ with $(1,y_2,y_2+y_3)$.
The vector fields $Y_1 = \frac{\partial}{\partial y_1}$,
$Y_2 = \frac{\partial}{\partial y_2} + y_1 \frac{\partial}{\partial y_3}$
and $Y_3 = \frac{\partial}{\partial y_3}$ are well defined in $G/H$
and satisfy $[Y_1, Y_2] = Y_3$, $[Y_1, Y_3] = [Y_2, Y_3] = 0$.
Set $Z = \frac{\partial}{\partial z}$: the vector fields
$Y_1, Y_2, Y_3, Z$ form a basis for the tangent space at any point of
$M = (G/H) \times (-\epsilon, \epsilon)$.

\subsection{Proof of theorem \ref{theo:L-unstable}}

Our action $\theta$ is given by a one parameter family of matrices $A(z)$
with $A(0) = I$, $A^\sharp(0)$ nilpotent, where
\[ A^\sharp(z) = \begin{pmatrix} a_{11}'(z) & a_{12}'(z) \\
a_{21}'(z) & a_{22}'(z) \end{pmatrix}. \]
Set $B = A'(0)$: notice that $b_{13} = b_{23} = 0$ and
$b_{33} = a_{11}'(0) a_{22}(0) + a_{11}(0) a_{22}'(0)
- a_{12}'(0) a_{21}(0) - a_{12}(0) a_{21}'(0) = a_{11}'(0) + a_{22}'(0) =
\trace A^\sharp(0) = 0$.
Write $\hat A(z) = I + zB$:
clearly, $\hat a_{33}(z) =
\hat a_{11}(z) \hat a_{22}(z) - \hat a_{12}(z) \hat a_{21}(z)$
for all $z$ and this family of matrices therefore defines
a family of automorphisms $\hat\phi_z: G \to G$ and an action
\begin{align} \hat\theta: G \times M &\to M \notag\\
(g, (g_1H, z)) &\mapsto (\hat\phi_z(g) g_1 H, z) \notag \end{align}
for which
\begin{align}
\hat X_1(y_1,y_2,y_3,z) &=
(1+zb_{11}) Y_1 + zb_{21} Y_2 + zb_{31} Y_3 \notag\\
\hat X_2(y_1,y_2,y_3,z) &=
zb_{12} Y_1 + (1+zb_{22}) Y_2 + zb_{32} Y_3 \notag\\
\hat X_3(y_1,y_2,y_3,z) &= Y_3 \notag 
\end{align}
Since $\theta$ and $\hat\theta$ are $C^1$ close in a neighborhood of $z = 0$,
it suffices to construct perturbations of $\hat\theta$.

Let $(p_2, -p_1)$ be a nonzero vector in the kernel of $A^\sharp(0)$, i.e.,
\begin{align}
b_{11} p_2 - b_{12} p_1 &= 0, \notag\\
b_{21} p_2 - b_{22} p_1 &= 0. \notag
\end{align}
Let $(q_1, q_2, q_4) \in \RR^3$ be a solution of the equation
\[ q_1 p_1 + q_2 p_2 + q_4 (b_{32} p_1 - b_{31} p_2) = 0 \]
with $q_4 \ne 0$ (which is possible since $p_1$ and $p_2$ are not both zero)
and define $c_{ij} = q_i p_j$.
Set
\begin{align}
\tilde X_1(y_1,y_2,y_3,z) &= \hat X_1(y_1,y_2,y_3,z) +
\psi(z) \left( c_{11} Y_1 + c_{21} Y_2 + c_{41} Z \right), \notag\\
\tilde X_2(y_1,y_2,y_3,z) &= \hat X_2(y_1,y_2,y_3,z) +
\psi(z) \left( c_{12} Y_1 + c_{22} Y_2 + c_{42} Z \right), \notag\\
\tilde X_3(y_1,y_2,y_3,z) &= \hat X_3(y_1,y_2,y_3,z), \notag
\end{align}
where $\psi$ is a bump function with $\psi(0) > 0$.
A straightforward computation verifies that 
\[ [\tilde X_1, \tilde X_2] = \tilde X_3, \quad
[\tilde X_1, \tilde X_3] = [\tilde X_2, \tilde X_3] = 0, \]
so that this defines an action.
Since either $c_{41} \ne 0$ or $c_{42} \ne 0$,
orbits of this action intersecting $z = 0$ are not compact.
This is the required perturbation.

\subsection{Proof of theorem \ref{theo:mixed} (instability part)}

As above, take an explicit action $\hat\theta$ with correct $A'(0)$:
\begin{align}
\hat X_1(y_1,y_2,y_3,z) &= e^{-\lambda z} Y_1, \notag\\
\hat X_2(y_1,y_2,y_3,z) &= e^{2\lambda z} Y_2, \notag\\
\hat X_3(y_1,y_2,y_3,z) &= e^{\lambda z} Y_3. \notag
\end{align}
A straightforward computation verifies that, for
\begin{align}
\tilde X_1(y_1,y_2,y_3,z) &= \hat X_1(y_1,y_2,y_3,z), \notag\\
\tilde X_2(y_1,y_2,y_3,z) &= \hat X_2(y_1,y_2,y_3,z)
- \frac{c}{\lambda} e^{\lambda z} Z, \notag\\
\tilde X_3(y_1,y_2,y_3,z) &= \hat X_3(y_1,y_2,y_3,z) - c Y_1, \notag
\end{align}
we have
\[ [\tilde X_1, \tilde X_2] = \tilde X_3, \quad
[\tilde X_1, \tilde X_3] = [\tilde X_2, \tilde X_3] = 0 \]
and this therefore defines a local action $\tilde\theta$.
For any nonzero $c$, the orbits of $\tilde\theta$
crossing $z=0$ are noncompact.
This, for small $c$, is the required perturbation.

\section{Unique ergodicity of flows \\ in a nilpotent manifold}

Let $X = v_1 Y_1 + v_2 Y_2 + v_3 Y_3$ with $v_2/v_1 \notin \QQ$:
we prove that the vector field $X$ is uniquely ergodic in $G/H$.

In order to prove this result we consider the diffeomorphism
$\phi: \TT^2 \to \TT^2$ defined by $\phi(z,w) = (\alpha z, zw)$
where $\alpha = e^{2\pi c}$, $c \notin \QQ$
(here we identify $\TT^2 = \Ss^1 \times \Ss^1$, $\Ss^1 \subset \CC$):
this is also uniquely ergodic.
We employ the ideas of Van der Corput's theorem (\cite{KN}).

\subsection{A diffeomorphism of the torus}

\begin{prop}
\label{prop:torus}
The diffeomorphism $\phi: \TT^2 \to \TT^2$
defined by $\phi(z,w) = (\alpha z, zw)$
where $\alpha = e^{2\pi c}$, $c \notin \QQ$, is uniquely ergodic.
\end{prop}

Notice that the usual unit measure $\mu$ (the Haar measure,
a constant multiple of the Lebesgue measure) is invariant under $\phi$.
This proposition implies that all orbits of $\phi$ are dense in $\TT^2$.

\proof
We must prove (\cite{Mane} or Theorem 9.2 of \cite{Mane2}) that if
$f: \TT^2 \to \CC$ is continuous then the sequence
\[ g_N = \frac{1}{N}(f + f \circ \phi + \cdots + f \circ \phi^{N-1}) \]
converges uniformly to a constant (the integral of $f$).
Furthermore, it suffices to consider $f$ a Laurent polynomial in $z$ and $w$.
Indeed, let $f$ be an arbitrary continuous function and let $f_n$ be a sequence
of Laurent polynomial tending uniformly to $f$.
If $a$ is the integral of $f$ and $a_n$ is the integral of $f_n$
we clearly have $a_n \to a$.
We may subtract $a$ from $f$ and $a_n$ from $f_n$ and assume without
loss of generality that $a = a_n = 0$.
Given $\epsilon > 0$, let $n_0$ be such that
$\| f_{n_0} - f \|_{\infty} < \epsilon/2$.
We thus have
\[ \frac{1}{N} \left\|
{
(f_{n_0} + f_{n_0} \circ \phi + \cdots + f_{n_0} \circ \phi^{N-1})
- (f + f \circ \phi + \cdots + f \circ \phi^{N-1})
}\right\|_{\infty} < \frac{\epsilon}{2} \]
for any $N$.
Since $f_{n_0}$ is a Laurent polynomial, we assume that 
\[ \frac{1}{N} \left\|{
(f_{n_0} + f_{n_0} \circ \phi + \cdots + f_{n_0} \circ \phi^{N-1})
}\right\|_{\infty} < \frac{\epsilon}{2} \]
for sufficiently large $N$.
This implies 
\[ \|g_N\|_\infty = \frac{1}{N} \left\|{
(f + f \circ \phi + \cdots + f \circ \phi^{N-1})
}\right\|_{\infty} < \epsilon \]
for large $N$, as desired.

Since a Laurent polynomial is a sum of monomials,
it is sufficient to prove our claim for $f(z,w) = z^n w^m$.
It is easy to show by induction that
\[ f \circ \phi^j(z,w) = \alpha^{jn + \frac{1}{2}j(j-1)m} z^{n+jm} w^m \]
for all $j$.
Our claim is trivial for $n = m = 0$.
For $m = 0$ and $n \ne 0$,
\begin{align}
g_N(z,w)
&= \frac{1}{N}
(z^n + \alpha^nz^n + \alpha^{2n}z^n + \cdots + \alpha^{(N-1)n}z^n) \notag\\
&= \frac{1}{N} \frac{\alpha^{Nn} - 1}{\alpha^n - 1} z^n.
\notag\end{align}
Since $\alpha$ and $n$ are fixed and $|z| = |\alpha| = 1$,
$(\alpha^{Nn} - 1)/(\alpha^n - 1) z^n$ is bounded
and $g_N$ tends uniformly to zero as $N$ goes to infinity.

Finally, if $m \ne 0$ we have
\begin{align}
g_N(z,w) &= \frac{1}{N}
(z^n + \alpha^nz^{n+m} + \cdots + \alpha^{jn+\frac{1}{2}j(j-1)m}z^{n+jm} + \notag\\
&\phantom{=} \cdots + \alpha^{(N-1)n + \frac{1}{2}(N-1)(N-2)m}z^{n+(N-1)m}) w^m.
\notag\end{align}
Let $\ww$ and $\zz$ be the unit vectors in $\CC^N$ given by
\begin{align}
\ww &= \frac{1}{\sqrt{N}}(1,\ldots,1), \notag\\
\zz &= \frac{1}{\sqrt{N}}(z^n,\ldots,\alpha^{jn+\frac{1}{2}j(j-1)m}z^{n+jm},
\ldots,\alpha^{(N-1)n + \frac{1}{2}(N-1)(N-2)m}z^{n+(N-1)m}).
\notag\end{align}
With the usual inner product, $|g_N(z,w)| = |\langle \zz, \ww \rangle|$
and we must prove that $\langle \zz, \ww \rangle$ tends to zero
uniformly in $z$.
Let $R_N: \CC^N \to \CC^N$ be the linear transformation
with $R_N(a_0, a_1, \ldots, a_{N-2}, a_{N-1}) =
(a_1, \ldots, a_{N-2}, a_{N-1}, a_0)$;
we have $R_N^\ast(a_0, a_1, \ldots, a_{N-2}, a_{N-1}) =
(a_{N-1}, a_0, a_1, \ldots, a_{N-2})$
and $R_N^\ast \ww = \ww$ so that
$\langle R_N^k \zz, \ww \rangle = \langle \zz, \ww \rangle$.
A straightforward computation yields
\begin{align}
\langle R^k_N \zz, \zz \rangle &= \frac{1}{N}
\Big( \frac{\alpha^{(N-k)km} - 1}{\alpha^{km} - 1}
\alpha^{kn + \frac{1}{2} k(k-1)m} z^{km} + \notag\\
&\phantom{=} +
(\alpha^{-\frac{1}{2}(N-k)(N-k-1)m} + \cdots + \alpha^{-\frac{1}{2}(N-1)(N-2)m})
\alpha^{(N-k)(-n)} z^{(N-k)(-m)}\Big); 
\notag\end{align}
for fixed $k$, the expression between parenthesis on the right hand side
is uniformly bounded in $z$ and therefore,
for any fixed $k$, $\langle R^k_N \zz, \zz \rangle$
tends to zero uniformly in $z$.

Let $\epsilon > 0$ and $M > 4/\epsilon^2$.
From what we just proved, there exists $N_0 > M$ such that
\[ |\langle R_N\zz, \zz \rangle|, |\langle R^2_N\zz, \zz \rangle|, \ldots
|\langle R^{M-1}_N\zz, \zz \rangle| < \frac{1}{M} \]
for all $N > N_0$. Then
\begin{align}
|\zz + R_N \zz + \cdots + R^{M-1}_N \zz| &=
\sqrt{\langle (\zz + \cdots + R^{M-1}_N \zz),
(\zz + \cdots + R^{M-1}_N \zz) \rangle} \notag\\
&\le \sqrt{M + \frac{1}{M}(M^2 - M)} \le \sqrt{2M}.
\notag\end{align}
Thus
\begin{align}
|\langle \zz, \ww \rangle| &= \frac{1}{M}
|\langle \zz + R_N \zz + \cdots + R^{M-1}_N \zz, \ww \rangle| \notag\\
&\le \frac{1}{M} |\zz + R_N \zz + \cdots + R^{M-1}_N \zz| |\ww| \notag\\
&\le \frac{1}{M} \sqrt{2M} = \sqrt{\frac{2}{M}} < \epsilon 
\notag\end{align}
for all $N > N_0$ and for all $z$.
\qed

\subsection{Vector fields on $G/H$}

\begin{theo}
\label{theo:ergodicflux}
Let $X$ be the vector field on $G/H$,
\[ X = v_1 Y_1 + v_2 Y_2 + v_3 Y_3 \]
where $v_1, v_2, v_3$ are constant real numbers.
If $v_2/v_1$ is irrational then $X$ is uniquely ergodic.
\end{theo}

Notice that this implies that every orbit of $X$ is dense in $G/H$.

\proof
Assume without loss of generality that $v_1 = 1$.
Let $\phi: \RR \times G/H \to G/H$
be the flux corresponding to $X$:
\[ \phi(t, (y_1,y_2,y_3)) =
(y_1 + t, y_2 + v_2 t, y_3 + v_3 t + v_2 y_1 t + \frac{1}{2} v_2 t^2). \]
Let $m$ be the Lebesgue or Haar measure on $G/H$:
it is clear that $X$ preserves $m$.
Given a Borel set $A$ and $x \in G/H$, let
\[ \tau(x,A) = \lim_{T \to \infty}
\frac{1}{T} \lambda(\{t \in [0, T] | \phi(t,x) \in A\}) \]
where $\lambda$ is the Lebesgue measure on $\RR$.
In order to prove that $X$ is uniquely ergodic
it suffices to prove that $\tau(x,A) = m(A)$
for every $x \in G/H$ and every $A$ in a basis
of the $\sigma$-algebra of Borel sets.
Our basis consists of sets of the form
\[ A = A_{\epsilon, y_1, D} =
\{ \phi(s,(y_1, y_2, y_3)), s \in (0,\epsilon), (y_2, y_3) \in D \} \]
where $\epsilon \in (0,1)$ is a real number
and $D \subset \TT^2$ is a Borel set (possibly a disc):
$m(A_{\epsilon, y_1, D}) = \epsilon \mu(D)$
where $\mu$ is the unit Lebesgue measure on $\TT^2$.

Let $t_0$ be the smallest nonnegative real number for which
$\phi(t_0,x)$ is of the form $(y_1, \ast, \ast)$:
set $\phi(t_0,x) = (y_1, \tilde x)$, $\tilde x \in \TT^2$.
Clearly,
\[ \tau(x, A_{\epsilon, y_1, D}) =
\lim_{T \to \infty} \frac{\epsilon}{T}
\card\{s \in \ZZ, 0 \le s < T, \phi_1^s(\tilde x) \in D \} \]
where $\phi_1: \TT^2 \to \TT^2$ is defined by $\phi$;
more precisely, 
\[ \phi(1,(y_1,y_2,y_3)) = (y_1,\phi_1(y_2,y_3)) =
(y_1, y_2 + v_2, y_3 - y_2 + v_3 + (y_1 - 1/2)v_2). \]
The diffeomorphism of the torus $\phi_1$ is conjugate
to a diffeomorphism of the form discussed in \ref{prop:torus}
and therefore uniquely ergodic.
Thus $\tau(x,A_{\epsilon,y_1,D}) = \epsilon \mu(D) = m(A)$
and we are done.
\qed

\section{Foliations}

\subsection{Excellent foliations and \\
near-identity nilpotent pseudogroups}

We define in this section excellent foliations of
$M = G/H \times (-\epsilon,\epsilon)$,
which have noncompact orbits but are near compact foliations.
We then define the horizontal subalgebra at each point;
we use results from \cite{BS}.
The reader should also see \cite{FF} for the more natural problem
of studying nilpotent groups of diffeomorphisms of one-manifolds.

Let $\Ff_0$ be the foliation of $M$
whose leaves are of the form $G/H \times \{\ast\}$.
Thus, $\Ff_0$ is the foliation defined by a horizontal homogeneous action.
We will study foliations $\Ff$ near $\Ff_0$.

We call $\Ff$ {\it acceptable} if its leaves are always transversal
to vertical lines, i.e., lines of the form
$\{\ast\} \times (-\epsilon, \epsilon)$.
Clearly, foliations near $\Ff_0$ are acceptable.
Alternatively, we say that a foliation is acceptable if the restriction
of the projection $\pi_M: M \to G/H$ to a leaf is a local diffeomorphism.

Given a function $\gamma: [a, b] \to G$ with $\gamma(a) = g_0$ and
a point $(g_0H, z_0) \in M$, we call the lift of $\gamma$ to be 
a path of the form $\tilde\gamma: [a, a+\delta) \to M$ such that
$\tilde\gamma(a) = (g_0H, z_0)$, $\tilde\gamma'(t)$ is tangent to $\Ff$
for all $t$ and $\pi_M \circ \tilde\gamma = \pi_H \circ \gamma$,
where $\pi_H: G \to G/H$ is the canonical projection.
Clearly, lifts exist for sufficiently small $\delta$
and there exists a unique maximal lift;
we say that $\gamma$ {\it lifts} if the domain
of the maximal lift $\tilde\gamma$ is $[a,b]$.
The foliation $\Ff$ is said to be {\it $k$-good} if
$\gamma$ lifs provided $|z_0|<\epsilon/2$ and the length of the image
of $\gamma$ is smaller than $k$;
abusing notation, we say that $\Ff$ is good if it is $k$-good
for some implicitly given value of $k$.

For a good foliation $\Ff$ and $g_0 \in G$, define the functions
$f_i^\pm: (-\epsilon/2, \epsilon/2) \to (-\epsilon, \epsilon)$
as follows. First, for $i = 1, 2, 3$, $s = \pm 1$, construct paths
$\gamma_{g_0,i}^s: [0,1] \to G$ defined by
\[ \gamma_{g_0,1}^s(t) = g_0 \begin{pmatrix} 1 & & \\
st & 1 & \\ 0 & 0 & 1 \end{pmatrix},
\gamma_{g_0,2}^s(t) = g_0 \begin{pmatrix} 1 & & \\
0 & 1 & \\ 0 & st & 1 \end{pmatrix},
\gamma_{g_0,3}^s(t) = g_0 \begin{pmatrix} 1 & & \\
0 & 1 & \\ st & 0 & 1 \end{pmatrix}. \]
Now, for $z \in (-\epsilon/2, \epsilon/2)$,
lift each $\gamma_{g_0,i}^s$ starting from $(g_0H, z)$ to obtain
points $(g_0H, f_{g_0,i}^s(z))$.
Clearly, each $f_{g_0,i}^s$ is smooth and strictly increasing.
Notice that $f_{g_0,i}^{-1}$ is the inverse of $f_{g_0,i}^{1}$;
from now on, we write $f_{g_0,i}$ instead of $f_{g_0,i}^1$.
Moreover, if $g_0H = g_1H$ then
$f_{g_1,3} = f_{g_0,3}$ and, for $i = 1, 2$,
\[ f_{g_1,i} = f_{g_0,i} \circ (f_{g_0,3})^{n_i} \]
for integers $n_1, n_2$:
this means that the choice of $g_0$ is not strongly relevant
and from now on $g_0$ will be omitted in the notation.
Finally,
\[ [f_2, f_1] = f_2^{-1} \circ f_1^{-1} \circ f_2 \circ f_1 = f_3,
[f_3, f_1] = [f_3, f_2] = id. \]

The functions $f_1, f_2, f_3$ form
a {\it near-identity nilpotent pseudogroup}
of functions as defined in \cite{BS}
if they are sufficiently near the identity;
in this case, we call $\Ff$ {\it excellent}.
A leaf $\Ll$ is called {\it abelian} if $f_3(z) = z$ for $(g_0H,z) \in \Ll$;
the foliation $\Ff$ is called {\it abelian} if all its leaves are abelian,
or, equivalently, if $f_3$ is the identity.

Let $\alpha: G \to \RR$ be a group homomorphism:
the $1$-form $d\alpha$ can be lifted to $G/H$
and therefore to $G/H \times \RR$,
where it is closed but probably not exact.
The $1$-form $dz$, where $z$ is the $\RR$-coordinate,
is exact in $G/H \times \RR$ and therefore $d\alpha + dz$ is a closed
$1$-form in $G/H \times \RR$, defining a foliation $\Ff_\alpha$.
We say that the homomorphism $\alpha$ is {\it rational}
if it takes elements of $H$ to rational numbers
and $\alpha$ is {\it irrational} otherwise.
The theorem below is a slightly modified version of theorem 2 from \cite{BS}.

\begin{theo}
\label{theo:foliation}
Given a manifold of the form $G/H$ there exists $\epsilon_{G/H}$ such that,
if $\Ff_1$ is a smooth foliation of $G/H \times (-1,1)$
with $d_{C^1}(\Ff_0,\Ff_1) < \epsilon_{G/H}$,
then one of the three conditions hold:
\begin{enumerate}
\item{$\Ff_1$ has at least one compact leaf and is an abelian foliation;
for any maximal connected open set $U$
of $G/H \times (-1,1)$ containing no compact leaf,
there exists a homomorphism $\alpha: G \to \RR$,
an open set $\Jj \subseteq G/H \times \RR$
and a homeomorphism $\Phi: \Jj \to U$
taking $\Ff_\alpha$ to $\Ff_1$.}
\item{$\Ff_1$ has no compact leaf and is an abelian foliation; 
there exists an irrational homomorphism $\alpha: G \to \RR$,
an open set $\Jj \subseteq G/H \times \RR$
and a homeomorphism $\Phi: \Jj \to G/H \times (-1,1)$
taking $\Ff_\alpha$ to $\Ff_1$.}
\item{$\Ff_1$ has no compact leaf and has an abelian leaf $\Ll_0$;
there exists a rational homomorphism $\alpha: G \to \RR$,
an open set $\Jj \subseteq G/H \times \RR$
and a homeomorphism $\Phi: \Jj \to G/H \times (-1,1)$
taking one leaf $\hat\Ll_0$ of $\Ff_\alpha$ to a leaf $\Ll_0$.}
\end{enumerate}
\end{theo}

The construction of the homeomorphism $\Phi$ is such that the $G/H$ coordinate
is preserved, i.e., $\Phi(gH,\ast) = (gH,\ast)$.

\subsection{The horizontal subalgebra}

For a group homomorphism $\alpha: G \to \RR$,
let $\alpha': \Gg \to \RR$ be the induced Lie algebra homomorphism.
We call the kernel $\Ee \subset \Gg$ of $\alpha'$
the {\it horizontal subalgebra}.
This can also be contructed as follows:
let $a$ be the translation number of $f_2$ relative to $f_1$;
the horizontal subalgebra $\Ee \subset \Gg$ is defined as the set of matrices
\[ \begin{pmatrix} 0 & 0 & 0 \\ -ay & 0 & 0 \\ z & y & 0 \end{pmatrix}. \]
In order to present another interpretation for $\Ee$
we first construct a local action $\bF: G \times M \to M$.

For $p_0 = (gH, z_0) \in M$ and $h \in G$,
let $\delta: [0,1] \to G$ be a continuous path with
$\delta(0) = I$, $\delta(1) = h$.
Construct a path $\gamma: [0,1] \to M$
whose image is contained in a leaf of $\Ff$
(or, equivalently, such that $\gamma'(t)$ is tangent to $\Ff$ for all $t$)
such that $\gamma(0) = p_0$, $\gamma(t) = (\delta(t)gH, \ast)$;
define $\bF(h,p_0) = \gamma(1)$.
It is easy to prove independence on the choice of $\delta$
and that $\bF$ is indeed a local action, i.e., that
$\bF(h_1h_2,p_0) = \bF(h_1,\bF(h_2,p_0))$.
Let $\Psi$ be the derivative of $\bF$, i.e., $\Psi_{p_0}(v) = \bF(\exp(v),p_0)$.
It follows easily from the theorem above that if $p_0$ is
sufficiently far from the boundary of $M$
then $\Psi_{p_0}(v)$ is defined for all $v \in \Ee_1$,
where $\Ee_1$ is the closed neighborhood
of radius $1$ of $\Ee$; by compactness, for such $p_0$, the function
$\Psi_{p_0}: \Ee_1 \to M$ is uniformly continuous.

From now on we shall consider two cases,
depending on whether $\alpha$ is a {\it rational} or {\it irrational}
homomorphism, or, equivalently, depending on whether $a$,
the translation number of $f_2$ relative to $f_1$,
is rational or irrational.
The irrational case corresponds either to case 2
in theorem \ref{theo:foliation} or to case 1 with $\alpha$ irrational.
In the irrational case, the closure of $\Psi_{p_0}(\Ee)$ is
$M_{1,p_0} = \Phi(\ast,\hat z_0)$ where $\Phi^{-1}(p_0) = (gH,\hat z_0)$,
a topological manifold homeomorphic to $G/H$.
The vector field $-a Y_1 + Y_2$ is tangent to $M_{1,p_0}$ and,
from theorem \ref{theo:ergodicflux}, uniquely ergodic in $M_{1,p_0}$.
The topological manifold $M_{1,p_0}$ thus receives a unit measure $\mu$.

Consider the $(y_1,y_2)$ torus $\TT^2 = \RR^2/\ZZ^2$,
where $\RR^2$ has coordinates $y_1$ and $y_2$.
In the rational case, the projection of $\Psi_{p_0}(\Ee) \subset \RR^2$
onto $\TT^2$
is a rational curve and the projection onto $G/H$ is a $2$-dimensional torus.
If $p_0$ belongs to an abelian leaf (for instance, in case 1),
then $\Psi_{p_0}(\Ee)$ is a $2$-dimensional torus $M_{2,p_0}$ contained
in $\Phi(\ast,\hat z_0)$ where $\Phi^{-1}(p_0) = (gH,\hat z_0)$.
If $p_0$ belongs to a nonabelian leaf, 
however, there may exist infinitely many points $(gH,z)$
in $\Psi_{p_0}(\Ee)$ above each point $gH$ of the torus.
For each $gH$, let $z^+$ and $z^-$ be the supremum and infimum, respectively, 
of the set of values of $z$ for which $(gH,z)$
belongs to $\Psi_{p_0}(\Ee)$.
The tori $M_{2,p_0}^+$ and $M_{2,p_0}^-$ are formed by points of the form
$(gH,z^+)$ and $(gH,z^-)$, respectively, and are contained in
abelian noncompact orbits.
We shall usually take $p_0$ as an abelian leaf to avoid these complications
but in either case
the vector fields $Y_3$ and $-a Y_1 + Y_2$ are tangent to $M_{2,p_0}$ or $M_{2,p_0}^\pm$
and define unique unit invariant measures $\mu$ or $\mu^\pm$ in
$M_{2,p_0}$ or $M_{2,p_0}^\pm$.

If $p_0 = (gH,z_0)$ and $p_1 = (gH,z_1)$ belong to the same leaf
then projection yields a measure preserving homeomorphism
$L_{p_0,p_1}: M_{\ast,p_0} \to M_{\ast,p_1}$.

\vfil
\section{Stability}

\subsection{The function $\tau_\Ab$}

Let $\Ab: G \to \RR^2$ and $\Ab': \Gg \to \RR^2$ be the abelianization
of the group $G$ and its Lie algebra $\Gg$, respectively, i.e.,
the projections defined by
\[ \Ab\left( \begin{pmatrix} 1 & 0 & 0 \\ y_1 & 1 & 0 \\
y_3 & y_2 & 1 \end{pmatrix} \right) = (y_1, y_2), \quad
\Ab'\left( \begin{pmatrix} 0 & 0 & 0 \\ v_1 & 0 & 0 \\
v_3 & v_2 & 0 \end{pmatrix} \right) = (v_1, v_2). \]

Let $\theta$ the {\em homogeneous horizontal} action
\[\begin{matrix}
\theta: & G \times M & \to & M \\
& (g, (g_1H, z)) & \mapsto & (\phi_z(g) g_1H, z)
\end{matrix}\]
where $M = G/H \times (-\epsilon, \epsilon)$ and
$\phi_z$ is a smooth family of automorphisms of $G$
such that $\phi_0$ is the identity.
The orbits of $\theta$ form the horizontal foliation $\Ff_0$.
Let $\tilde\theta$ be a smooth local action near $\theta$;
its orbits form a foliation $\Ff$ of $M$.
If $\tilde\theta$ is a perturbation of $\theta$,
we may assume $\Ff$ to be an excellent foliation
and use the results of the previous section.
For a vector $v \in \Gg$, write $\tilde\tau_{p_0}(v)$
for the vector $w$ in $\Gg$ taken by $\tilde\theta$ to 
a vector of the form $(v,\ast)$:
\[ \bT(t) = \tilde\theta(\exp(tw),p_0), \quad \bT'(0) = (v,\ast). \] 
Given $p_0$, define $\tau_{\Ab,p_0}: \Ee \to \RR^2$ by
\[ \tau_{\Ab,p_0}(v) = \int_{M_{\ast,p_0}} \Ab'(\tilde\tau_p(v)) d\mu(p) \]
where $M_{\ast,p_0} = M_{1,p_0}$, $M_{\ast,p_0} = M_{2,p_0}$
or $M_{\ast,p_0} = M^+_{2,p_0}$ depending on whether
the translation number is irrational or rational.

In order to provide a geometric interpretation for the function $\tau_{\Ab}$,
we introduce an auxiliary function $\xi_{p_0}: \Ee_1 \to G$:
\[ \tilde\theta(\xi_{p_0}(v),p_0) = \Psi_{p_0}(v); \]
$\xi_{p_0}$ is uniformly continuous in $\Ee_1$.

\begin{lemma}
\label{lemma:tauxi}
For a noncompact perturbation $\tilde\theta$ of $\theta$,
\[ \tau_{\Ab,p_0}(v) = \lim_{t \to \infty} \frac{1}{t} \Ab(\xi_{p_0}(tv)). \]
\end{lemma}

\proof
We first consider the irrational case. Take 
\[ v = (v_1,v_2,v_3) =
\begin{pmatrix} 0 & 0 & 0 \\ v_1 & 0 & 0 \\ v_3 & v_2 & 0 \end{pmatrix}
\in \Ee; \]
we have $\Ab(\xi_{p_0}(tv)) = \int_0^t \Ab'(\tilde\tau_{p(s)}(v)) ds$
where $p(s) = \Psi_{p_0}(sv)$.
If $(v_1, v_2) \ne 0$ it follows from theorem \ref{theo:ergodicflux}
and from Birkhoff's theorem that 
\[ \lim_{t \to \infty} \frac{1}{t} \Ab(\xi_{p_0}(tv)) =
\int_{M_{1,p_0}} \Ab'(\tilde\tau_p(v)) d\mu(p) = \tau_{\Ab,p_0}(v), \]
as required.
If $v = (0,0,1)$, we have
\[ \lim_{t \to \infty} \frac{1}{t} \Ab(\xi_{p_0}(tv)) =
\int_{S(p_0)} \Ab'(\tilde\tau_p(v)) dp, \]
where $S(p_0) \subset M_1$ is the circle
$\bF(\exp(t (0,0,1)),p_0)$.
We have to show that if $p_1 \in M_1$ then
\[ \int_{S(p_0)} \Ab'(\tilde\tau_p(v)) dp =
\int_{S(p_1)} \Ab'(\tilde\tau_p(v)) dp, \]
i.e., that
\[ \lim_{t \to \infty} \frac{1}{t} \Ab(\xi_{p_0}(tv)) =
\lim_{t \to \infty} \frac{1}{t} \Ab(\xi_{p_0}(tv+w)) \]
where $v \in \Ee$ and $w \in \Ee_1$.
Since we have already seen that both limits exist,
it sufficies to prove that $\Ab(\xi_{p_0}(tv)) - \Ab(\xi_{p_0}(tv+w))$ is bounded,
and this follows from uniform continuity of $\xi_{p_0}$.

In the rational case, $M^+_{2,p_0}$ is a torus.
If the orbit of $v$ is dense,
the result follows directly from Birkhoff's theorem, as before.
Otherwise, we again have that the limit equals an integral over a circle
and independence from the choice of circle follows from uniform continuity
of $\xi_{p_0}$.
\qed

\subsection{Dependence of $\tau_{\Ab,p_0}$ on $p_0$}

Set 
\[ E_3 = \begin{pmatrix} 0 & 0 & 0 \\ 0 & 0 & 0 \\ 1 & 0 & 0 \end{pmatrix}. \]

\begin{lemma}
\label{lemma:tauab3}
The value of $\tau_{\Ab,p}(E_3)$ is constant
on each abelian orbit of $\tilde\theta$.
If $\tilde\theta$ has a compact orbit then
$\tau_{\Ab,p}(E_3) = 0$ for any point $p$.
\end{lemma}

\proof
We first prove that $\Ab(\xi_{p}(E_3))$ is constant on abelian orbits.
Indeed, let $p_0 = (g_0H,z_0)$ and $p_1 = (g_0H,z_1)$
belong to the same abelian orbit of $\tilde\theta$
with $\tilde\theta(g,p_0) = p_1$.
The closed loops $\tilde\gamma_{p_0,3}$ and $\tilde\gamma_{p_1,3}$
are the images under $\tilde\theta$ of paths $\gamma_0, \gamma_1 \in G$
with $g^{-1} \gamma_1(t) g = \gamma_0(t)$.
Clearly, $\xi_{p_i}(E_3) = \gamma_i(1)$ at $p_i$, $i = 0, 1$
and therefore $\Ab(\xi_{p_0}(E_3)) = \Ab(\xi_{p_1}(E_3))$.
From the proof of lemma \ref{lemma:tauxi}, it follows
that this constant value of $\Ab(\xi_{p}(E_3))$ is $\tau_{\Ab,p_0}(E_3)$.

At a compact orbit, the path $\tilde\gamma_{p_0,3}$
is the commutator of $\tilde\gamma_{p_0,1}$ and $\tilde\gamma_{p_0,2}$
and therefore $\xi_{p_0}(E_3)$ is also a commutator and therefore
$\Ab(\xi_{p_0}(E_3)) = \tau_{\Ab,p_0}(E_3) = 0$.
On the other hand, if $p_0$ belongs to a noncompact abelian orbit
but a compact orbit exists then there exist $p_1, p_2, \ldots$
is the same orbit as $p_0$ tending to $p_\infty$, a point in a compact orbit.
We have $\Ab(\xi_{p_i}(E_3)) = \Ab(\xi_{p_0}(E_3))$ for all $i$.
On the other hand, the closed loops $\tilde\gamma_{p_i,3}$
tend uniformly to $\tilde\gamma_{p_\infty,3}$ and therefore
$\lim_{i \to \infty} \Ab(\xi_{p_i}(E_3)) = \Ab(\xi_{p_\infty}(E_3)) = 0$
and we have $\tau_{\Ab,p_0}(E_3) = \Ab(\xi_{p_0}(E_3)) = 0$.
\qed

The example presented in subsection 3.2 shows that $\tau_{\Ab,p_0}(E_3)$
is not always zero, not even if the foliation is abelian.
In that example, it is easy to verify that 
\[ \xi_{(gH,z)}(E_3) = \begin{pmatrix} 1 & 0 & 0 \\
c & 1 & 0 \\ e^{-\lambda z} & 0 & 1 \end{pmatrix} \]
and therefore $\tau_{\Ab,(gH,z)}(E_3) = (c,0) \ne 0$.

\begin{lemma}
\label{lemma:tauab}
If $p_0 = (g_0H,z_0)$ and $p_1 = (g_0H,z_1)$
belong to the same abelian orbit of $\tilde\theta$
then $\tau_{\Ab,p_1}(v) - \tau_{\Ab,p_0}(v)$
is a multiple of $\tau_{\Ab,p}(E_3)$.
\end{lemma}

\proof
Set $p_1 = \bF(h,p_0)$ and $F_h$ be
the local diffeomorphism $F_h(p) = \bF(h,p)$.
For $u \in \Ee$, write $u^h = h^{-1} u h$, the conjugate of $u$ by $h$.
We clearly have $h \exp(t u^h) = \exp(t u) h$ and therefore
$F_h(\Psi_{p_0}(t u^h)) = \Psi_{p_1}(t u)$.
From the definition of $\xi$, we have
$F_h(\tilde\theta(\xi_{p_0}(t u^h),p_0)) =
\tilde\theta(\xi_{p_1}(t u),p_1)$ and therefore
$|\xi_{p_0}(t u^h) - \xi_{p_1}(t u)|$ is bounded when $t$ goes to infinity,
yielding $\tau_{\Ab,p_0}(u^h) = \tau_{\Ab,p_1}(u)$.
Since $u^h - u$ is a multiple of $E_3$ we are done.
\qed

The following two results are now immediate.

\begin{coro}
\label{coro:tauabcomp}
If $\tilde\theta$ has a compact orbit then $\tau_{\Ab,p}$
is constant in each connected component of the complement of compact orbits.
\end{coro}

\begin{coro}
\label{coro:tauabimage}
The image of the map $\tau_{\Ab,p}$ has dimension $1$ and is constant
in each connected component of a vertical line of the form $(g_0H,\ast)$
minus the intersection with compact orbits.
\end{coro}

\subsection{Proof of theorems \ref{theo:L-stable} and \ref{theo:T-stable}}

We begin with theorem \ref{theo:L-stable}.
Assume by contradiction that there exists a non-compact action
$\tilde\theta$ near $\theta$ coinciding with $\theta$ outside
a neighborhood of an orbit.
Take $p_0 = (g_0H,z_0)$ and $p_1 = (g_0H,z_1)$, $z_1 > z_0$,
in the same non-compact orbit.
From corollary \ref{coro:tauabcomp}, $\tau_{\Ab,p_0} = \tau_{\Ab,p_1}$.
Let $v \in \Ee$ with $|\Ab'(v)| = 1$.
By definition of $\tau_{\Ab,p}$, we have
\[ \tau_{\Ab,p_i}(v) = \int_{M_{\ast,p_i}} \Ab'(\tilde\tau_p(v)) d\mu(p) \]
and therefore
\[ \int_{M_{\ast,p_0}}
(\Ab'(\tilde\tau_{L_{p_0,p_1}(p)}(v)) - \Ab'(\tilde\tau_{p}(v))) d\mu(p) = 0. \]
On the other hand, for each $p \in M_{\ast,p_0}$,
\[ \Ab'(\tilde\tau_{L_{p_0,p_1}(p)}(v)) - \Ab'(\tilde\tau_{p}(v)) =
\Ab' \left( \int_{\tilde z_0}^{\tilde z_1} 
\frac{d}{dz} \tilde\tau_{(gH,z)} dz \right) (v). \]
Notice that for the original action $\theta$, 
\[ \Ab' \tilde\tau_{(gH,z)}(v) = B(z) (\Ab'(v)),
\qquad
B(z) = \begin{pmatrix} a_{11}(z) & a_{12}(z) \\
a_{21}(z) & a_{22}(z) \end{pmatrix}^{-1} \]
and therefore
\[ \Ab' \frac{d}{dz} \tilde\tau_{(gH,z)}(v) \approx
- (B(z) A^\sharp(z) B(z))(\Ab'(v)). \]
Since $A^\sharp(z)$ is assumed to be invertible,
$- (B(z) A^\sharp(z) B(z))(\Ab'(v)) \ne 0$;
by compactness, for $\tilde\theta$ sufficiently near $\theta$,
$\Ab' \frac{d}{dz} \tilde\tau_{(gH,z)}(v)$ belongs to a convex open set
not containing the origin and therefore the integral above belongs
to the open cone spanned by this open set, contradicting the fact
that the integral equals zero.
This completes the proof of theorem \ref{theo:L-stable}.

For the proof of theorem \ref{theo:T-stable},
assume without loss of generality that $A^\sharp(u) \wedge u > 0$
for all nonzero vectors $u$.
We know that if $\tau_{\Ab,p_0}(v) = w$
then $\tau_{\Ab,p_1}(v) = cw$ for some positive number $c$.
We have $w \wedge \tau_{\Ab,p_1}(v) = 0$
and therefore
\[ w \wedge \Ab' \left( \int_{\tilde z_0}^{\tilde z_1} 
\frac{d}{dz} \tilde\tau_{(gH,z)} dz \right) (v) = 0. \]
For $\tilde\theta$ sufficiently near $\theta$,
we may assume that $\Ab'(\tilde\tau_p(v))$ always belongs
to a narrow angle around $w$;
since $B(z)$ is near the identity we have
$w \wedge (- (B(z) A^\sharp(z) B(z))(\Ab'(v)) ) > 0$,
again yielding a contradiction.

{\obeylines
\parskip=0pt
\parindent=0pt
Tania M. Begazo
Universidade Federal Fluminense
Instituto de Matem\'atica
Niterói, RJ 24020-140, Brazil
tania@mat.uff.br

\smallskip

Nicolau C. Saldanha
Depto. de Matem{\'a}tica, PUC-Rio
R. Mq. de S. Vicente 225
Rio de Janeiro, RJ 22453-900, Brazil
nicolau@mat.puc-rio.br
}

\end{document}